\documentclass[12pt]{amsart}
\usepackage{amssymb,amsthm}

\def\silentfootnote#1{{\let\thefootnote\relax\footnotetext{#1}}}

\theoremstyle{plain}
\newtheorem{thm}[equation]{Theorem}
\newtheorem{lem}[equation]{Lemma}
\newtheorem{prop}[equation]{Proposition}
\newtheorem{cor}[equation]{Corollary}
\numberwithin{equation}{section}

\theoremstyle{remark}
\newtheorem{rem}[equation]{Remark}

\theoremstyle{definition}

\DeclareMathOperator{\Hom}{Hom}

\DeclareMathOperator{\ann}{ann}

\DeclareMathOperator{\End}{End}
\DeclareMathOperator{\Ker}{Ker}

\DeclareMathOperator{\soc}{-soc}

\DeclareMathOperator{\op}{op}
\DeclareMathOperator{\st}{st}

\newcommand{\ba}{\begin{array}}
\newcommand{\ea}{\end{array}}

\title{Clifford correspondence for algebras}
\author{Sarah J.\ Witherspoon}
\address{Department of Mathematics and Statistics, University of
Massachusetts, Amherst, MA 01003}
\address{(2001--02)  
Department of Mathematics and Computer Science, Amherst College,
Amherst, MA 01002}
\email{wither@math.umass.edu}
\date{July 20, 2001}

\begin{document}

\silentfootnote{Research supported in part by National Security Agency
Grant \#MDA904-01-1-0067.}

\begin{abstract}
We give a Clifford correspondence for an algebra $A$ over an algebraically
closed field, that is an algorithm for constructing some
finite-dimensional simple $A$-modules
from simple modules for a subalgebra and endomorphism algebras.  
This applies to all finite-dimensional simple $A$-modules in the case
that $A$ is finite-dimensional and semisimple with a given semisimple
subalgebra.
We discuss connections between our work and earlier results, with a view
towards applications particularly to finite-dimensional semisimple
Hopf algebras.
\end{abstract}

\maketitle

\section{Introduction}
\label{sec:intro}

Clifford theory is an important collection of results, due to A.\ H.\
Clifford \cite{clifford37}, relating
representations of finite groups to representations of normal subgroups.
One of the results, the Clifford correspondence, is a one-to-one
correspondence between simple modules for a group and pairs
$(V,W)$ where $V$ is a simple module for a given normal subgroup $H$ 
of finite index, and $W$
is a simple module for a twisted group algebra of a subgroup of
$G/H$.  (These pairs are taken up to a standard equivalence.)  The
correspondence is constructive, and may be used to build simple
modules of $G$.

Many authors have generalized this result to other types of
groups and rings, for example to locally compact groups 
\cite{mackey58, rieffel79a}, and to group-graded rings \cite{dade70, fell69,
ward68}.  Others have extended at least part of the theory to 
arbitrary rings \cite{R}, non-normal subgroups \cite{A}, and Hopf algebras
\cite{S2, VZ, W}.

The present work is an attempt to unify some of these generalizations,
and to fill in some gaps, with the goal of finding future applications.  
We extend some of the results of 
Rieffel for arbitrary rings \cite{R}, at the expense of assuming the 
rings are algebras over algebraically closed fields.  
Statements and proofs are presented module-theoretically,
providing an alternative view of Rieffel's theory.
We connect our results to
Alperin's version of the Clifford correspondence for non-normal subgroups
\cite{A}, which in fact also holds for arbitrary rings.

This note applies in particular to Hopf algebras, extending some of
the work in Clifford theory for Hopf algebras.  This may be the
most promising potential application of our work:  In the case
of finite-dimensional semisimple Hopf algebras, the full theory is valid,
and the constructions may be concrete enough to classify all simple modules in
some cases.
There is still much to be learned about finite-dimensional semisimple Hopf
algebras $H$.  For example a conjecture of Kaplansky, that the dimension
of a simple $H$-module should divide $\dim(H)$, has not been settled in this
full generality.

In \cite{MW} Clifford theory is used to prove Kaplansky's Conjecture
when $H$ is a finite-dimensional semisimple Hopf algebra constructed by a
sequence of crossed products involving group algebras and their duals;
in particular this includes all semisimple Hopf algebras of prime power
dimension in characteristic 0.
(Other partial proofs of this conjecture have been demonstrated by authors
using completely different methods \cite{etingof-gelaki98,nichols-richmond96}.)

As further motivation, we mention some other recent applications of the
Clifford correspondence.
In \cite{RR} Clifford theory is used to show that the simple modules for
many Hecke algebras and groups of interest may all be obtained from those of
the affine Hecke algebras of type A.
Clifford theory arises in \cite{BW} in a classification of the
finite-dimensional simple modules for
Hopf algebra extensions of some of the down-up algebras
(these Hopf algebras are essentially generalizations of the Borel subalgebra of
$U_q(\mathfrak{sl}_3)$).

Let $k$ be an algebraically closed field.  All algebras will be algebras
over $k$, and a subalgebra will have the same identity element as the algebra.
All modules will be left modules unless otherwise specified.

\section{Static modules and Alperin's correspondence}

In this section we summarize some results from \cite{A} that will be needed.
Let $A$ be a ring with identity, $M$ a fixed left $A$-module, and
$E=\End_A(M)^{\op}$.  
Our functions will act on the left, so that $M$ is a {\it right} module for
$E$.
(Alperin chooses to have his functions act on the right, and
accordingly takes $E$ to be the endomorphism algebra itself.)
If $N$ is another $A$-module, then the space $\Hom_A(M,N)$
is naturally a (left) $E$-module under composition of functions (as $E$ 
takes the opposite multiplication from an endomorphism algebra).  We say that
$N$ is {\it static} (with respect to $M$) if the $A$-homomorphism 
\begin{equation}\label{phi1}
\phi:M\otimes_E \Hom_A(M,N)\rightarrow N
\end{equation}
given by $\phi(m\otimes f) =f(m)$ is an isomorphism. 
A left $E$-module $U$ is {\it static} (with respect to $M$) if the
$E$-homomorphism 
\begin{equation}\label{psi}
\psi: U\rightarrow \Hom_A(M,M\otimes_E U)
\end{equation}
given by $\psi(u)(m)=m\otimes u$ is an isomorphism.  
As $M\otimes_E -$ and $\Hom_A(M, -)$ are adjoint functors, there is a
category equivalence between the full subcategories of static $A$-modules
and static $E$-modules:

\begin{prop}[{\it [1, Lemma 1]}]\label{a}
The category of static $A$-modules is equivalent with the category of
static $E$-modules.  The equivalence is given by the functors $\Hom_A(M, -)$
and $M\otimes_E -$.
\end{prop}

Alperin uses this category equivalence to give new proofs of Auslander 
equivalence and Morita equivalence, as well as to
give a version of the Dade-Cline equivalence for non-normal subgroups
of a finite group \cite{A}.
(The Dade-Cline equivalence \cite{C,D} is a category equivalence which
implies the Clifford correspondence for group-graded rings.)
In Section 3, we will 
use Proposition \ref{a} to extend Rieffel's results further for algebras
over $k$.

Next we give Alperin's version of a Clifford correspondence for 
non-normal subgroups, valid for arbitrary rings, which will be used
in Section 4.
If $M$ is a fixed $A$-module,
we say that another $A$-module $N$ has an $M$-{\it presentation} if there is an
exact sequence
$$ \oplus ' M\rightarrow \oplus M\rightarrow N\rightarrow 0
$$
of $A$-modules.
Letting $B$ be a subring of $A$,
we say this presentation is $B$-{\it split} if it splits on restriction
to $B$.

Now let $V$ be a $B$-module and $M$ the $A$-module induced from $V$, that is
$M=V\!\uparrow^A=A\otimes_B V$ where $A$ acts by multiplication on the first
factor.  Let $E=\End_A(V\!\uparrow^A)^{\op}$ as above, and 
$F=\End_B(V\!\uparrow^A)^{\op}$, so that $E\subseteq F$.  If $U$ is a left 
$E$-module, we define the {\it restriction} of $U$ to 
$F$ to be the $F$-module $F\otimes_E U$ (as in \cite{A}), 
denoted $U\!\downarrow_F$.

\begin{prop}[{\it [1, Theorem]}] \label{non-normal}
The category of $A$-modules with $B$-split $V\!\uparrow^A$-presentations 
is equivalent
with the category of $E$-modules whose restrictions to $F$ are projective.
The equivalence is given by the functors $\Hom _A(V\!\uparrow^A,-)$ and
$V\!\uparrow^A\otimes_E-$.
\end{prop}

Alperin's proof of the proposition is valid in this general setting, and
it uses Proposition \ref{a}.  The proof involves verifying that:

\begin{itemize}
\item[(i)] If $N$ is an $A$-module with a $B$-split 
$V\!\uparrow^A$-presentation, then
$N$ is static (with respect to $V\!\uparrow^A$) and 
$\Hom_A(V\!\uparrow^A,N)\!\downarrow_F$ is
$F$-projective.

\item[(ii)] If $U$ is an $E$-module whose restriction to $F$ is projective,
then $V\!\uparrow^A\otimes_EU$ has a $B$-split $V\!\uparrow^A$-presentation, 
and $U$ is static (with respect to $V\!\uparrow^A$).
\end{itemize}
These assertions can be verified exactly as in \cite{A}.

\section{Stable Clifford correspondence}

There are traditionally two parts to the Clifford correspondence.
One part involves the case where a module for a normal subgroup is stable
in the whole group, and the other part involves the connection between
modules for a stabilizer subgroup and modules for the whole group.
Here we will deal with the stable part for algebras more generally;
the rest appears in the next section.

Let $B\subset A$ be algebras over the algebraically closed field $k$, and
$V$ a $B$-module.  
We say that $V$ is $A$-{\it stable} if the restriction of $V\!\uparrow^A
=A\otimes_BV$ to $B$, denoted $V\!\uparrow^A\downarrow_B$, is
isomorphic to a finite direct sum of copies of $V$.

\begin{lem}\label{b}
Let $V$ be a finite-dimensional simple $A$-stable $B$-module, 
$M=V\!\uparrow^A$, and $E=\End_A(M)^{\op}$.
\begin{itemize}
\item[(i)] A finite-dimensional $A$-module $N$ is static with respect to $M$
if and only if
$N\!\downarrow_B$ is isomorphic to a direct sum of copies of $V$.

\item[(ii)] All finite-dimensional $E$-modules are static with respect to $M$.
\end{itemize}
\end{lem}

\begin{proof}
(i) First we claim that if $N\!\downarrow_B\cong V^{\oplus n}$, then the
$B$-homomorphism
$$\gamma: V\otimes_k \Hom_A(M,N)\rightarrow N\!\downarrow_B$$
given by $\gamma(v\otimes f)=f(1\otimes v)$ is an isomorphism, where 
$B$ acts only on the left factor of $V\otimes_k \Hom_A(M,N)$.
To see this, note that $\Hom_A(M,N)\cong \Hom_B(V,N\!\downarrow_B)$ 
as induction and restriction are adjoint functors.  
Since $N\!\!\downarrow_B\cong V^{\oplus n}$, the latter is isomorphic to
$k^{\oplus n}$ by
Schur's Lemma, from which it is clear that $\gamma$ is surjective.
This also shows that the dimensions of the modules are the same, so
$\gamma$ is bijective.

In particular, we may let $N=M$ as $V$ is $A$-stable, and we obtain 
\begin{equation}\label{vem}
  V\otimes_kE\cong M
\end{equation}
as left $B$-modules and right $E$-modules (where $E$ acts only on the right
factor of $V\otimes_k E$).
Therefore if $N$ is any finite-dimensional $A$-module with $N\!\downarrow_B
\cong V^{\oplus n}$, we have $B$-module isomorphisms
\begin{equation}\label{eqn}
M\otimes_E\Hom_A(M,N)\stackrel{\sim}{\rightarrow} V\otimes_k\Hom_A(M,N)
\stackrel{\sim}{\rightarrow} N\!\downarrow_B.
\end{equation}
It is straightforward to check that the composite map is precisely the map
$\phi$ of (\ref{phi1}), an $A$-homomorphism.

Conversely, if $N$ is static then by (\ref{vem}),
$$N\!\downarrow_B\cong (M\otimes_E\Hom_A(M,N))\!\downarrow_B\cong 
V\otimes_k\Hom_A(M,N)$$
as $B$-modules.
As $B$ acts only on the left factor in $V\otimes_k\Hom_A(M,N)$, this is
isomorphic to a direct sum of copies of $V$.

(ii) Let $U$ be a finite-dimensional $E$-module. Then by (\ref{vem}), we have
\begin{eqnarray*}
\Hom_A(M,M\otimes_EU) &\cong & \Hom_B(V, (M\otimes_EU)\!\downarrow_B)\\
   &\cong & \Hom_B(V,(V\otimes_k U)\!\downarrow_B)\\
  &\cong & \Hom_B(V,V^{\oplus \dim U}) \cong  k^{\oplus \dim U}.
\end{eqnarray*}
It follows that the map $\psi:U\rightarrow \Hom_A(M,M\otimes_EU)$ of
(\ref{psi}) is surjective.
As these modules also have the same dimension, $\psi$ must be an isomorphism.
Therefore $U$ is static. \end{proof}

Combining the lemma with Proposition \ref{a} yields:

\begin{thm}\label{stable-a}
Let $V$ be a finite-dimensional simple $A$-stable $B$-module and 
$E=\End_A(V\!\uparrow^A)^{\op}$.
The category of finite-dimensional $E$-modules is equivalent with the
category of finite-dimensional $A$-modules whose restriction to
$B$ is a direct sum of copies of $V$.
The equivalence is given by the functors 
$V\!\uparrow^A\otimes_E -$ and $\Hom_A(V\!\uparrow^A, -)$. 
\end{thm}

We point out that if $U$ is a finite-dimensional $E$-module, 
then 
\begin{equation}\label{module}
  V\!\uparrow^A\otimes_EU\cong V\otimes_k U
\end{equation}
as $B$-modules (where $B$ acts on the left factor only) by (\ref{vem}).
The $B$-module $V\otimes_kU$ may then be given an $A$-module structure via
this isomorphism. 
This may give a simpler description of such $A$-modules, and
it potentially allows an easy count of their dimensions.

We next show that any finite-dimensional {\it simple} $A$-module containing a
simple $A$-stable $B$-module $V$ on restriction to $B$ may be described in
this way:

\begin{lem}\label{cc}
Let $N$ be a finite-dimensional simple $A$-module such that $N\!\downarrow_B$
contains a copy of a finite-dimensional simple $A$-stable $B$-module $V$.
Then $N\!\downarrow_B$ is a direct sum of copies of $V$.
\end{lem}

\begin{proof}
Define $\pi:V\!\uparrow^A\rightarrow N$ by $\pi(a\otimes v)=a\cdot v$ 
(where we have identified $V$ with a $B$-submodule of $N$).
Then $\pi$ is a nonzero $A$-homomorphism, and so must be surjective, as
$N$ is simple.
Therefore $N\!\downarrow_B$ is a homomorphic image of 
$V\!\uparrow^A\downarrow_B$, a direct sum of copies of $V$,
and so is itself such a direct sum.
\end{proof}

In light of Lemma \ref{cc}, we immediately have a corollary
of Theorem \ref{stable-a}:

\begin{cor}[Stable Clifford correspondence]\label{stable}
Let $V$ be a finite-dimensional simple $A$-stable $B$-module, and 
$E=\End_A(V\!\uparrow^A)^{\op}$.
There is a one-to-one correspondence between finite-dimensional simple
$E$-modules and finite-dimensional simple
$A$-modules containing $V$ on restriction to $B$.
\end{cor}

\begin{rem} This generalizes the stable Clifford correspondence of finite
group theory:  Let $G$ be a finite group with normal subgroup $H$.  
Let $A=kG$ and $B=kH$.  A $B$-module $V$ is $A$-stable precisely when it is
$G$-{\it stable}, that is all its conjugates under $G$ are isomorphic.
(If $g\in G$, the {\it conjugate} $B$-module ${}^gV$ has underlying vector
space $V$, and module structure given by $h\cdot_gv=(g^{-1}hg)\cdot v$ for
all $h\in H,\ v\in V$.)
In this case, $E$ is isomorphic to a twisted group algebra of $G/N$.
(In general, $V\!\uparrow^A\downarrow_B$ will be isomorphic to a direct sum
of copies of conjugates of $V$, one for each element of $G/N$, and $E$ will
be isomorphic to a twisted group algebra of a subgroup of $G/N$.)
\end{rem}

We now relate these results to those of Rieffel in \cite{R}.
We say that an ideal $J$ of $B$ is $A$-{\it invariant} if $AJ=JA$.
We say that $B$ is {\it normal} in $A$ if for every
(two-sided) ideal $I$ of $A$, $I\cap B$ is $A$-invariant.
(We caution that Rieffel's definition of normal is somewhat different
in case $A$ and $B$ are not both semisimple Artinian.)

\begin{lem}\label{invariant}
Let $B\subset A$ be normal, $V$ a finite-dimensional simple
$B$-module, $J=\ann_B V$, and assume $V\!\uparrow^A$ is finite-dimensional.
Then $V$ is $A$-stable if and only if $J$ is $A$-invariant.
\end{lem}

\begin{proof}
Let $M=V\!\uparrow^A$, $I=\ann_AM$, and assume $V$ is $A$-stable.
We claim that $J=I\cap B$.
Note first that $I\cap B\subseteq J$ as $V\subseteq M$.
Since $M\!\downarrow_B$ is a direct sum of copies of $V$,
$J$ annihilates $M$, so $J\subseteq I\cap B$.
Therefore $J=I\cap B$.
As $B$ is normal in $A$, it follows that $J$ is $A$-invariant.

Assuming now that $J$ is $A$-invariant, we again claim that $J=I\cap B$.
We have $I\cap B\subseteq J$ as $V\subseteq M$.
Now
$$ J\cdot M =JA\otimes_BV=AJ\otimes_BV=0$$
as $J$ is $A$-invariant and annihilates $V$.  
Therefore $J\subseteq I\cap B$, and we have $J=I\cap B$.
It follows that $B/(I\cap B) =B/J\cong \End_k(V)$ has a unique
simple module up to isomorphism, namely $V$.
As $I$ annihilates $M$, $M\!\downarrow_B$ is naturally a $B/(I\cap B)$-module,
and so must be a direct sum of copies of $V$.
\end{proof}

\begin{rem} Rieffel gives a result analogous to Corollary 
\ref{stable}, under the hypotheses that $A$ and $B$ are semisimple
Artinian rings and $B$ is normal in $A$ \cite[Proposition 2.15]{R}.
He uses the hypothesis that $J$ is $A$-invariant where we require
$V$ to be $A$-stable.
As applied to algebras over $k$, Corollary \ref{stable} is
more general. 
\end{rem}

\section{General Clifford correspondence}

In this section we discuss the case in which a finite
dimensional simple $B$-module $V$ is not necessarily $A$-stable.
We will assume that $A$ and $B$ are {\it finite-dimensional
semisimple} algebras over $k$, as the results and proofs are more 
transparent in this case.  

We will give a module-theoretic definition of stabilizer $S$ for $V$
and two new proofs of the full Clifford correspondence in this setting.
(Rieffel has a more ring-theoretic definition of stabilizer, and a full
Clifford correspondence for semisimple Artinian rings \cite[\S 2]{R}.)
The first of our proofs is traditional, involving one step to get simple
$S$-modules, and another to get from simple $S$-modules to simple
$A$-modules. The second proof uses Alperin's correspondence (Proposition
\ref{non-normal}), and bypasses stabilizers altogether.

Let $B$ be a normal subring of $A$ and $V$ a finite-dimensional
simple $B$-module.
A {\it stabilizer} (or {\it stability} subring)
for $V$ is a semisimple algebra
$S$ such that $B\subseteq S\subseteq A$, $\ B$ is a normal subring of
$S$, $V$ is $S$-stable, and $V\!\soc (V\!\uparrow^A)=V\!\soc(V\!\uparrow^S)$.
(Here, the $V$-{\it socle} of a module is the sum of all its $B$-submodules 
that are isomorphic to $V$.)
This definition is equivalent to that of Rieffel \cite[Definition 2.11]{R}
under the hypotheses that
$A$ and $B$ are semisimple, as the following lemma shows.
The lemma will be used in the first proof of the full Clifford correspondence.

\begin{lem}\label{stabilizer}
Assume $A$ and $B$ are finite-dimensional semisimple algebras with
$B\subset A$ normal. Let $V$ be a finite-dimensional simple 
$B$-module and $J=\ann_BV=pB$ where $p$ is a primitive central idempotent
of $B$.
Then a semisimple subalgebra $S$ of $A$, containing $B$ as a normal
subring, is a stabilizer for $V$ if and only if $J$ is $S$-invariant
and $S+AJ+JA=A$.
\end{lem}

\begin{proof}
Let $S$ be a stabilizer for $V$. By Lemma \ref{invariant}, $J$ is 
$S$-invariant.  
Note that $B/J\cong V^{\oplus \dim V}$ as $B$-modules, so that
$$
  V\!\soc (A\otimes _B(B/J))\cong V\!\soc(S\otimes_B(B/J))=S\otimes _B(B/J)
$$
as $S$ is a stabilizer for $V$.  Also note that as a $B$-module,
$(1-p)A\otimes_B(B/J)$ is a direct sum of copies of $V$, since $p$
annihilates it.
This means that the $B$-submodule $(1-p)A\otimes_B(B/J)$ of
$A\otimes_B(B/J)$ is contained in the image of the canonical inclusion
$$
  S\otimes_B(B/J)\hookrightarrow A\otimes_B(B/J).
$$
Applying the isomorphisms $S\otimes_B(B/J)\cong S/SJ$ (given by
$s\otimes(b+J)\mapsto sb+SJ$) and
$A\otimes_B(B/J)\cong A/AJ$, we now have that $(1-p)A/AJ$ embeds in
$S/SJ$ under the canonical map, that is $(1-p)A\subseteq S+AJ$.  Therefore
$$
  A=(1-p)A +pA \subseteq S+AJ+JA,
$$
so that $S+AJ+JA=A$ as desired.

Conversely, suppose that $S$ is a semisimple subalgebra 
of $A$ containing $B$ as a
normal subring, $J$ is $S$-invariant, and $S+AJ+JA=A$.
Then $V$ is $S$-stable by Lemma \ref{invariant}.  As $J$ annihilates $V$,
we have 
\begin{eqnarray*}
  V\!\uparrow^A\downarrow_B = (A\otimes_BV)\!\downarrow_B & = & ((S+AJ+JA)
\otimes_BV)\!\downarrow_B\\
   &=& ((S+JA)\otimes_BV)\!\downarrow_B.
\end{eqnarray*}
Note that $(S\otimes_BV)\!\downarrow_B\subseteq
V\!\soc(V\!\uparrow^A)$
as $V$ is $S$-stable.  Further, as $J=pB$, $p$ acts as the identity on
$JA\otimes_BV$, which implies that $JA\otimes_BV$ intersects
$V\!\soc(V\!\uparrow^A)$ trivially.  Therefore
$$
  V\!\soc(V\!\uparrow^A)=S\otimes_BV=V\!\soc(V\!\uparrow^S).
$$
\end{proof}

\begin{rem} Rieffel shows that stabilizers exist (but are not unique)
in case both $A$ and $B$ are semisimple Artinian:
$S=B+(1-p)A(1-p)$ is a stabilizer for $V$ contained in every other
stabilizer, and $S=pAp+(1-p)A(1-p)$ (the centralizer of $p$)
is a stabilizer for $V$ containing every other stabilizer.
\end{rem}

\begin{rem} In case $H$ is a normal subgroup of a finite
group $G$, $A=kG$ and $B=kH$, the {\it stabilizer} subgroup (or {\it inertia}
subgroup) $L$ for a simple $B$-module $V$ consists of all
$g\in G$ such that ${}^gV\cong V$.
The algebra $S=kL$ is a stabilizer for $V$ in the above sense,
and is the unique subalgebra of $A$ arising from a subgroup of $G$
that is a stabilizer for $V$.
\end{rem}

We now have two expressions for the
endomorphism algebra $E=\End_A(V\!\uparrow^A)^{\op}$.  
For, as a $B$-homomorphism must map $V$ into the
$V$-socle of a module, we have vector space isomorphisms
\begin{equation}\label{endalg}
  \End_A(V\!\uparrow^A)^{\op} \cong \Hom_B(V,V\!\uparrow^A\!\downarrow_B)
   \cong  \Hom_B(V,V\!\uparrow^S\!\downarrow_B)
  \cong  \End_S(V\!\uparrow^S)^{\op}.
\end{equation}
It may be verified that the composite isomorphism preserves products in
the endomorphism algebras.

\begin{thm} \label{cc2}
Assume $A$ and $B$ are finite-dimensional semisimple algebras, with $B$
normal in $A$.  Let $V$ be a finite-dimensional simple $B$-module
and $S$ a stabilizer for $V$.  
There is a one-to-one correspondence between finite-dimensional simple
$S$-modules containing $V$ on restriction to $B$, and finite-dimensional 
simple $A$-modules containing $V$ on restriction to $B$.
The correspondence is given by induction of modules from $S$ to $A$.
\end{thm}

\begin{proof}
We may assume $S\neq A$.  First let $N$ be a finite-dimensional simple 
$S$-module with $V$ contained in $N\!\downarrow_B$.  
As $V$ is $S$-stable, $N\!\downarrow_B$ must be a direct sum of
copies of $V$ by Lemma \ref{cc}.
Let $M$ be any simple $A$-submodule
of $N\!\uparrow^A$.  Then $M$ is a direct summand of $N\!\uparrow^A$, and so
$$
  0\neq \Hom_A(N\!\uparrow^A,M)\cong\Hom_S(N,M\!\downarrow_S).
$$
As $N$ is a simple $S$-module, this implies that $N$ embeds in
$M\!\downarrow_S$.  On the other hand, if $J=\ann _BV=pB$, then by
Lemma \ref{stabilizer},
\begin{eqnarray*}
  N\!\uparrow^A\downarrow_S &=& ((S+AJ+JA)\otimes_SN)\!\downarrow_S\\
   &=& ((S+JA)\otimes_SN)\!\downarrow_S\\
  &=& N\oplus (JA\otimes_SN)\!\downarrow_S
\end{eqnarray*}
as $J$ (and so $p$) annihilates $N$ ($N\!\downarrow _B$ is a direct sum of 
copies of $V$), but $p$ acts as the identity on 
$JA\otimes_SN$.  This implies that there is a {\it unique} copy
of $N$ in $N\!\uparrow^A\!\downarrow_S$.  Now by semisimplicity we have
$N\!\uparrow^A \cong M\oplus W$
for some $A$-module $W$, but restricting to $S$, we see that $W$ cannot
contain $N$ (as $M\!\downarrow_S$ contains the unique copy of $N$ in
$N\!\uparrow^A\downarrow_S$).  Applying the above argument to any simple
$A$-submodule of $W$ leads to the contradiction that $W$ must also contain
$N$.  Therefore $M\cong N\!\uparrow^A$ and $N\!\uparrow^A$ is simple.

Now let $M$ be {\it any} simple $A$-module containing $V$ on restriction
to $B$.  Let $N$ be a simple $S$-submodule of the image of $S\otimes_BV$
under the map $\pi :A\otimes_B V\rightarrow M$ given by $\pi(a\otimes
v)=a\cdot v$. 
By the previous arguments, $N\!\uparrow^A$ is a simple $A$-module.
By the definition of $N$,
$$
 \Hom_A(N\!\uparrow^A, M)\cong \Hom_S(N,M\!\downarrow_S)\neq 0.
$$
As both $M$ and $N\!\uparrow^A$ are simple $A$-modules, this implies
$M\cong N\!\uparrow^A$.
\end{proof}

Now we are ready to give the full Clifford correspondence.

\begin{thm}[General Clifford correspondence]\label{mainthm}
Assume that $A$ and $B$ are finite-dimensional semisimple algebras with $B$
normal in $A$.
Let $V$ be a finite-dimensional simple $B$-module,
and $E=\End_A(V\!\uparrow^A)^{\op}$.  There is a one-to-one
correspondence between finite-dimensional simple $E$-modules and
finite-dimensional simple $A$-modules containing $V$ on restriction
to $B$.
\end{thm}

\begin{proof}
We give two proofs:

(i) Let $S$ be a stabilizer for $V$.
By Corollary \ref{stable} and (\ref{endalg}), there is a one-to-one
correspondence between finite-dimensional simple $E$-modules and 
finite-dimensional simple $S$-modules containing
$V$ on restriction to $B$.
Theorem \ref{cc2} provides the rest of the correspondence.

(ii) Let $F=\End_B(V\!\uparrow^A)^{\op}$.  As $V\!\uparrow^A\downarrow_B$ is
semisimple, $F$ is a direct sum of matrix algebras and so is 
semisimple itself.  Thus every $F$-module is projective.
Applying Proposition \ref{non-normal}, it remains to
show that each simple $A$-module $N$ containing $V$ on restriction to $B$
has a $V\!\uparrow^A$-presentation (automatically $B$-split as $B$ is
semisimple).  

Note that $N$ is a quotient of $V\!\uparrow^A$ via $\pi: V\!\uparrow^A
\rightarrow N$, $\pi(a\otimes v)=a\cdot v$ (where we have identified $V$
with a $B$-submodule of $N$).
Let $K$ be any simple component of the kernel of $\pi$, so that $K$ embeds
in $V\!\uparrow^A$.  By semisimplicity, $K$ is also a quotient of $V\!
\uparrow^A$.
Taking a direct sum of copies of $V\!\uparrow^A$, one for each simple
summand of the kernel of $\pi$, we obtain a $V\!\uparrow^A$-presentation
$$
  \oplus V\!\uparrow^A\rightarrow V\!\uparrow^A\rightarrow M\rightarrow 0.
$$
\end{proof}

Note that the second proof only requires $V\!\uparrow ^A$ and $V\!\uparrow
^A\!\downarrow_B$ to be semisimple $A$- and $B$-modules, respectively.
It does not require $A$ and $B$ to be semisimple algebras, nor does it
require $B$ to be normal in $A$.

The correspondence of the theorem
is given by the functors $\Hom_A(V\!\uparrow^A,-)$ and $V\!\uparrow^A\otimes
_E -$ from Proposition \ref{non-normal}.  
So the second proof {\it and} the construction of simple $A$-modules
(from $V$ and simple $E$-modules) bypasses a stabilizer $S$ altogether.
This may be useful in situations where stabilizers do not have
desired properties.  For example, in the case of Hopf algebras (discussed
in the next section), we do not know whether there is always a stabilizer 
that is a {\it Hopf} subalgebra.
On the other hand, in specific applications, it may be useful to know 
a stabilizer $S$
in order to obtain more detailed information about the modules.

\begin{rem} Rieffel gives a result analogous to Theorem
\ref{cc2}, but for {\it arbitrary} rings \cite[Theorem 3.25]{R}.
It is quite possible that our methods could be extended to include algebras
more generally, and in any case Theorem 3.25 of \cite{R} combined with
Corollary \ref{stable} should yield a version of Theorem \ref{mainthm}
valid under more general hypotheses.  We do not pursue this here as
the ring-theoretic definitions required to state Rieffel's Theorem 3.25
would lead us too far afield (see \cite[\S 3]{R}).
It would be interesting to find a module-theoretic statement and proof of
his theorem.
\end{rem}

\section{Clifford correspondence for Hopf algebras}

Let $B\subset A$ be a finite-dimensional $H$-Galois extension for a Hopf 
algebra $H$ over $k$ (see \cite{montgomery93} for definitions).
Let $V$ be a finite-dimensional $B$-module.
By \cite[Remark 3.2(4)]{S1}, $V$ is $A$-stable if and only if there 
is a left $B$-linear and right $H$-colinear isomorphism 
\begin{equation}\label{Phi}
  \Phi: A\otimes_BV\stackrel{\sim}{\rightarrow} V\otimes_kH.
\end{equation}
This is Schneider's definition of $A$-stable (which does not require any
finiteness assumptions, see \cite[\S 3]{S1}).
A stable Clifford correspondence in this context is given by Schneider
\cite[Remark 5.8(2)]{S2} and van Oystaeyen and Zhang \cite[Theorem 5.4]{VZ}. 
Their results are reformulated and extended in \cite{W} to a general
Clifford correspondence in case $H$ is cocommutative, using a stabilizer
Hopf subalgebra defined by Schneider in \cite{S1}.

In this section, we summarize some of the results of \cite{W} and show how 
they follow from this note.  
At the same time, our results here or those of Rieffel in \cite{R} give a
much more general version of the Clifford correspondence for Hopf algebras,
which may turn out to be useful.

Suppose that $A=B\#_{\sigma} H$ is a crossed product
algebra of a $k$-algebra $B$ with $H$ (a special case of an $H$-Galois
extension).
That is, $\sigma$ is a two-cocycle
mapping $H\otimes H$ to $B$ and there is 
a measuring $H\otimes B\rightarrow B$,
denoted $h\otimes b\mapsto h\cdot b$.
These maps must satisfy appropriate properties 
(see \cite[Chapter 7]{montgomery93}).
We first show that $B$ is a normal subring of $A$, as defined in Section 3. 

\begin{prop} The subalgebra $B$ of $A=B\#_{\sigma}H$ is a normal subring of
$A$.
\end{prop}

\begin{proof}  Let $I$ be any (two-sided) ideal of $A$.  We must show
that $A(I\cap B)=(I\cap B)A$.  Let $a=b\# h\in A$ and $i\in I\cap B$.
Then by definition of multiplication in $A$,
\begin{eqnarray*}  
  ai = (b\# h)(i\# 1) & = & \sum b(h_1\cdot i)\# h_2 \\
              &=& \sum (b(h_1\cdot i)\# 1)(1\# h_2).
\end{eqnarray*}
This is an element of $(I\cap B)A$, as by \cite[(7.2.4)]{montgomery93},
$h_1\cdot i\in I\cap B$.  Therefore $A(I\cap B)\subseteq (I\cap B)A$.
The other containment is proved similarly. 
\end{proof}

If $A$ is any finite-dimensional Hopf algebra with a {\it normal} Hopf 
subalgebra $B$ (see \cite[Definition 3.4.1]{montgomery93}),
then $A=B\#_{\sigma}H$ for some two-cocycle $\sigma$, where $H=A/AB^+$
($B^+=B\cap(\Ker \varepsilon)$) \cite{schneider92}, and so the lemma
applies. Alternatively, it can be shown directly that any normal Hopf
subalgebra of a Hopf algebra is a normal subring.

Under these hypotheses, the stable Clifford correspondence for
finite-dimensional Hopf algebras
\cite[Theorem 3.1 or Corollary 3.3]{W}
follows immediately from our Theorem \ref{stable-a} 
(or Corollary \ref{stable}).  
The proof of the Clifford correspondence
in \cite{W} is based on calculations involving the various 
structure maps, potentially allowing more detailed information
about the structure of $A$-modules (see \cite[Theorem 2.2(i)]{W}).  
The proof we give of Theorem \ref{stable-a} here is essentially a
generalization of those given for Hopf algebras in \cite{S2,VZ}, and is 
based instead on facts about the relevant categories.

Under the hypothesis that $H$ is {\it cocommutative}, we also
gave in \cite{W} a general Clifford correspondence.  
This is the case in which a stabilizer coalgebra $H_{\st}$ of a $B$-module
$V$, defined by Schneider, is a Hopf subalgebra of $H$.  We will review
this definition.

Let $V$ be a finite-dimensional $B$-module.  If $C\subseteq H$ is a 
subcoalgebra, let $A(C)=\Delta_A^{-1}(A\otimes C)$ where $\Delta_A:
A\rightarrow A\otimes H$ is the right $H$-comodule map arising from the
crossed product $A=B\#_{\sigma}H$.  We say that $C$ {\it stabilizes} $V$ if
$V\otimes C\cong A(C)\otimes_BV$ as left $B$-modules and right $C$-comodules
\cite[p.\ 216]{S1}.  Let $H_{\st}=\sum C$, the sum over all subcoalgebras
$C\subseteq H$ such that $C$ stabilizes $V$.  
Assuming $H$ is cocommutative, $H_{\st}$ is a Hopf
subalgebra of $H$ by \cite[Theorem 4.4]{S1}.  (Note that $H$ is pointed
as $k$ is algebraically closed and $H$ is cocommutative 
\cite[p.\ 76]{montgomery93}. In the more general case, $H_{\st}$ is a
subcoalgebra of $H$, but not necessarily a Hopf subalgebra.)

Let $S=A(H_{\st})=B\#_{\sigma}H_{\st}$.  By \cite[Remark 3.2(4), and Theorems 
4.4 and 5.4(1)(a)]{S1}, $S$ is a stabilizer for $V$ (our definition).
The cocommutative Clifford correspondence \cite[Corollary 3.5]{W} {\it in
the semisimple case} follows immediately from our Theorem \ref{mainthm},
and can be done in two steps (as in proof (i) of the theorem) involving
the stabilizer $S$.  The $A$-modules produced are of the form
$V\!\uparrow^A\otimes _E U\cong A\otimes_S(V\otimes_kU)$ for $E$-modules $U$
by (\ref{module}),
and as $S$ arises from a Hopf subalgebra of $H$, this may allow an explicit 
description of their structure (see for example \cite[Theorem 2.2(i)]{W}).
On the other hand, Theorem \ref{mainthm} gives a Clifford correspondence
more generally for semisimple $A$ and $B$ ($H$ need not be cocommutative),
regardless of whether any of the stabilizers $S$ arise from Hopf subalgebras
of $H$.

We make some final remarks that may be helpful in calculating examples.
If the measuring of $B$ by $H$ is trivial (that is $h\cdot b
=\varepsilon(h)b$ for all $h\in H$ and $b\in B$), then every $B$-module $V$
is $A$-stable.  Indeed, the map $\Phi$ of (\ref{Phi}) is given by $\Phi((b\# h)
\otimes _B v)=\rho(b)(v)\otimes h$ where $\rho: B\rightarrow \End_k(V)$
is the structure map of the $B$-module $V$.

Let $E=\End_A(V\!\uparrow^A)^{\op}$.  By \cite[Theorem 3.6]{S1} and Schur's
Lemma, $E$ is an $H$-crossed product over the field $k$.  That is, $E$ is
isomorphic to a twisted product $k_{\alpha}[H]=k\#_{\alpha}H$ as any
measuring of $k$ by $H$ is necessarily trivial. 
Therefore $E$-modules are equivalent to projective representations of
$H$ \cite[Proposition 2.4]{boca}.  
These may be studied via the methods of Boca in \cite{boca}.

\end{document}